\documentclass[10pt]{amsart}
\usepackage{amsmath}
\usepackage{graphicx}
\usepackage{latexsym}
\usepackage{color}
\usepackage{amscd}
\usepackage[all]{xy}
\usepackage{enumerate}
\usepackage{hyperref}
\usepackage{subfigure}
\usepackage{soul}
\usepackage{comment}
\usepackage{epsfig}
\usepackage{epstopdf}
\newtheorem{thm}{Theorem}

\newtheorem{theorem}[thm]{Theorem}

\newtheorem{corollary}[thm]{Corollary}

\theoremstyle{definition}

\newtheorem*{definition*}{Definition}

\newcommand{\CPb}{\overline{\mathbb{CP}}{}^{2}}
\newcommand{\CP}{{\mathbb{CP}}{}^{2}}

\newcommand{\K}{{\rm K3}}

\def \x {\times}
\def \eu{{\text{e}}}

\begin{document}

\title[Inequivalent Lefschetz fibrations on rational and ruled surfaces] 
{Inequivalent Lefschetz fibrations on \\ rational and ruled surfaces}

\author[R. \.{I}. Baykur]{R. \.{I}nan\c{c} Baykur}
\address{Department of Mathematics and Statistics, University of Massachusetts, Amherst, MA 01003-9305, USA}
\email{baykur@math.umass.edu}

\begin{abstract}
In this short note, we give an explicit construction of inequivalent  Lefschetz pencils and fibrations of same genera on blow-ups of all rational and ruled surfaces. This complements our earlier results in \cite{BaykurJSG}, concluding that \emph{every} symplectic $4$--manifold, after sufficiently many blow-ups, admits inequivalent Lefschetz pencils and fibrations, which cannot be obtained from one another even via any sequence of fibered Luttinger surgeries. 
\end{abstract}

\maketitle

\setcounter{secnumdepth}{2}
\setcounter{section}{0}

\section{Introduction} 

A pair of Lefschetz pencils or fibrations\footnote{We assume that a \emph{pencil}, unlike a Lefschetz fibration, always has base points. We also assume that a Lefschetz pencil or fibration always has critical points, and each singular fiber contains only one critical point and no exceptional spheres (i.e. it is \emph{relatively-minimal}).} on a closed oriented smooth \mbox{$4$--manifold} are said to be \emph{equivalent} (or \emph{isomorphic}) if there is an orientation-preserving self-diffeomorphism of the total space which commutes the two maps. Simon Donaldson's ground-breaking work in \cite{Donaldson} provides arbitrarily large genera \emph{pencils} on any symplectic \mbox{$4$--manifold.} A fundamental problem then is to determine if there are relatively-minimal inequivalent pencils or fibrations of \emph{the same genus and number of base points} on a given symplectic $4$--manifold. For brevity, we refer to them as \emph{inequivalent pencils} and \emph{inequivalent fibrations}, where it should be understood that their genera and number of base points are already the same. 

This note is an amendment to our results in \cite{BaykurJSG}, where we proved that any symplectic $4$-manifold, which is not a rational or a ruled surface, after a number of blow-ups, admits arbitrarily many inequivalent Lefschetz pencils and fibrations. Moreover, they are not related via \emph{fibered Luttinger surgeries} \cite{BaykurJSG}. Here we will give similar examples of pencils and fibrations on \emph{rational} and \emph{ruled surfaces}. Recall that rational and ruled surfaces are birationally equivalent to the complex projective plane and $S^2$--bundles over Riemann surfaces, so they are diffeomorphic to blow-ups of $\CP$, $S^2 \x \Sigma_h$, $S^2 \widetilde{\x} \Sigma_h$, which are known to have unique symplectic structures up to symplectomorphisms and deformations. 

Our examples in \cite{BaykurJSG} were obtained from Donaldson pencils via \emph{partial doubling sequences}, which involve blow-ups and degree doublings of pencils  introduced in \cite{Smith, AurouxKatzarkov}. These pencils were distinguished by looking at the collection of exceptional spheres, which arise as \emph{multisections} \cite{DonaldsonSmith, BaykurHayano} of different degrees in respective pencils. Thus, blowing-up all the base points, we obtain inequivalent Lefschetz fibrations. However, the same strategy does not work for rational or ruled surfaces, since in this case, one is not guaranteed to have disjoint multisections representing a given collection of  exceptional classes \cite{BaykurJSG}.

We will instead apply partial doublings to hand-picked pencils on rational and ruled surfaces in order to derive inequivalent pencils which have \emph{different numbers of reducible fibers} versus irreducible ones. The point here is that ``doublings are not created equal''; depending on which subcollections of base points we pick, we can arrive at pencils with topologically different singular fibers. Since the distinction comes from the number of reducible fibers, we get inequivalent fibrations after blowing-up all the base points. Our main result is:

\begin{theorem}
Any rational or ruled surface, after sufficiently many blow-ups, admits inequivalent relatively-minimal pencils and fibrations of the same genera and number of base points, which cannot be obtained from one another via any sequence of fibered Luttinger surgeries.
\end{theorem}

\noindent Last part of the theorem comes easy for our examples: performing a fibered Luttinger surgery is equivalent to conjugating a subword of the corresponding monodromy factorization of a pencil by some Dehn twist; see \cite{Auroux2, BaykurJSG}. This can be done whenever the twist curve is fixed, up to isotopy, by the mapping class given by the subword. Clearly, the number of reducible fibers is unchanged under this operation, so the examples with different numbers of reducible fibers we produce are also inequivalent up to fibered Luttinger surgeries. 

\smallskip
There are many other examples of inequivalent Lefschetz pencils and fibrations. A pair of inequivalent \emph{fibrations} on a blow-up of $T^2 \x \Sigma_2$, whose fibers have different divisibility in homology, was discussed by Ivan Smith in \cite{SmithThesis}. Several inequivalent fibrations on homotopy elliptic surfaces, distinguished by their monodromy groups,  were discovered by  Jongil Park and Ki-Heon Yun in \cite{ParkYun, ParkYun2}. As for inequivalent \emph{pencils}, first examples were constructed by Kenta Hayano and the author on \mbox{$4$--manifolds} homeomorphic to blow-ups of the $\K$ surface \cite{BaykurHayano}. The construction in \cite{BaykurHayano} uses pairs of lantern substitutions which modify pencils differently, but land on the same $4$--manifold. Notably, Noriyuki Hamada recently constructed pairs of inequivalent pencils on irrational ruled surfaces $S^2 \widetilde{\times} \Sigma_h$, which \emph{only} differ by the configuration of their base points on the reducible fibers, so they do not yield inequivalent fibrations after blow-ups. 


With our aforementioned result from \cite{BaykurJSG}, we arrive at the following conclusion, which demonstrates how diverse the Donaldson-Gompf correspondence between symplectic $4$--manifolds and Lefschetz pencils and fibrations is: 

\begin{corollary}
Any closed symplectic $4$--manifold, possibly after blow-ups, admits inequivalent relatively-minimal pencils and fibrations of the same genera and number of base points, which cannot be obtained from one another via any sequence of fibered Luttinger surgeries.
\end{corollary}

\smallskip

\section{Preliminaries} 

Here we will quickly review the basic notions and background results on Lefschetz pencils\,/\,fibrations, positive Dehn twist factorizations, fibered Luttinger surgeries and partial conjugations, and doublings of Lefschetz pencils. For further details on these topics, the reader can turn to \cite{GS} and \cite{ADK, AurouxKatzarkov, BaykurJSG}.

Let $X$ be a closed, oriented $4$-manifold, and $B=\{ b_j \}$, $C=\{ p_i \}$ be finite, non-empty sets of points in $X$.  A \emph{Lefschetz pencil} $(X,f)$ is given by a surjective map $f\colon X \setminus B \to S^2$, which is a submersion on $X \setminus C$, such that around each \emph{base point} $b_j$ and each \emph{critical point} $p_i$, it conforms to local complex models $(z_1,z_2) \mapsto z_1/z_2$ and $(z_1, z_2) \mapsto z_1 z_2$, respectively.  A \emph{Lefschetz fibration} is defined similarly for $B = \emptyset$. After blowing-up every base point $b_j$ of a pencil, we obtain a Lefschetz fibration with disjoint exceptional spheres $S_j$ as its sections. We say $(X,f)$ is a \emph{genus--$g$} pencil\,/\,fibration, for $g$ the genus of the regular fiber $F$. We will assume that every  fiber contains at most one critical point, which can always be achieved after a small perturbation. A \emph{singular fiber} has a nodal singularity at the critical point $p_i$, and is obtained by shrinking a simple loop $a_i$ on $F$, called the \emph{vanishing cycle}. We have a \emph{reducible} fiber (with two connected components) if $a_i$ is separating $F$, and \emph{irreducible} if $a_i$ is non-separating. 

Let $\Sigma_g^m$ denote a compact, oriented surface of genus $g$ with $m$ boundary components, where $\Sigma_g=\Sigma_g^0$. The \emph{mapping class group} $\Gamma_g^m$ is the group of orientation-preserving self-diffeomorphisms of $\Sigma_g^m$ which restrict to identity along $\partial \Sigma_g^m$, modulo isotopies of the same type. Let $t_a \in \Gamma_g^m$ denote the positive (right-handed) Dehn twist along a simple loop $a$ on $\Sigma_g^m$. A genus--$g$ Lefschetz pencil $(X,f)$ with $m$ base points prescribes a \emph{positive factorization} of the boundary multi-twist (or \emph{monodromy factorization}), given by a relation
\[ t_{\delta_1} \cdot \ldots \cdot t_{\delta_m} = t_{a_1} \cdot \ldots \cdot t_{a_\ell} \]
in $\Gamma_g^m$, where $\delta_1, \ldots, \delta_m$ are curves parallel to distinct boundary components of $\Sigma_g^m$ (which we will also denote the boundary component by),  $a_i$ is a vanishing cycle corresponding to the critical point $p_i$ of $f$, and $\ell=|C|$.  Conversely, given a positive factorization as above, one can construct such a Lefschetz pencil $(X,f)$. For $g \geq 2$, we indeed get a one-to-one correspondence between equivalence classes of  \mbox{Lefschetz pencils} and positive factorizations up to \textit{Hurwitz moves} (trading subwords $t_{a_i}t_{a_{i+1}}$ and $ t_{t_{a_i}(a_{i+1})}  t_{a_i}$) and \textit{global conjugations} (replacing every $t_{a_i}$ in the factorization with $t_{\phi{(a_i)}}$, for the same $\phi \in \Gamma_g^n$). All goes the same for Lefschetz \emph{fibrations} for $m=0$.

By the pioneering works of Donaldson and Gompf  \cite{Donaldson, GS}, every symplectic \mbox{$4$-manifold} admits a  Lefschetz pencil whose fibers are symplectic,  and conversely, the total space of a Lefschetz pencil or a fibration  (recall $C \neq \emptyset$) can always be equipped with a symplectic form for which all regular fibers are symplectic. 

As observed in \cite{ADK}, \emph{Luttinger surgery} (where one removes a neighborhood of a Lagrangian torus in a symplectic $4$--manifold and glues it back in differently so as to produce a new symplectic $4$--manifold) can be performed in a way compatible with the pencil structure,  when the restriction of the pencil map to the torus is a circle bundle over a loop on the base, away from the critical values \cite{Auroux2, ADK, BaykurJSG}. This surgery is equivalent to modifying the corresponding positive factorization by replacing a subword $t_{a_i} \cdot \ldots \cdot t_{a_k}$ with a conjugate $t_{\phi(a_i)} \cdot \ldots \cdot t_{\phi(a_k)}$ for some $\phi=t_{b}^{\pm}$, which can be done if and only if the mapping class $t_{a_i} \cdot \ldots \cdot t_{a_k}$ fixes $b$, up to isotopy \cite{Auroux2}. 

Another way to produce a new pencil from a given one, and on the same \mbox{$4$--manifold,} is the \emph{doubling} procedure, discussed for holomorphic pencils and Donaldson’s pencils in \cite{Smith}, and for pencils obtained via branched coverings of $\CP$ in \cite{AurouxKatzarkov}. In general, by \cite{BaykurJSG}[Lemma~3.1], one can double  \emph{any} (topological) genus--$g$ pencil $(X,f)$ with $m$ base points, provided $g \geq 2$ and $m \leq 2g-2$. The result is a genus--$g'$ pencil $(X, f')$ with $m'$ base points, where $g'=2g+m-1$ and $m'= 4m$. If the original pencil is compatible with a given symplectic form, then the resulting one is also compatible with a deformation equivalent symplectic form. 

\begin{figure}[ht!]
 \centering
     \includegraphics[width=10cm]{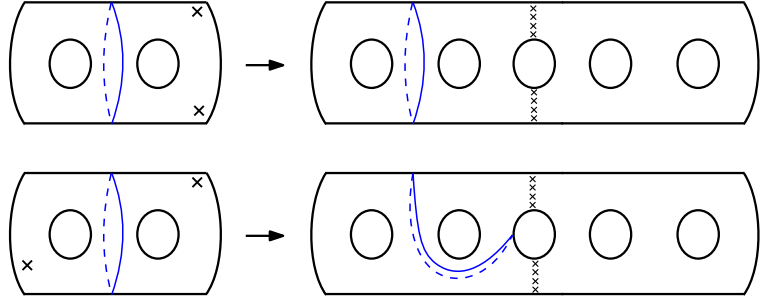}
     \caption{The lifts of a separating vanishing cycle (in blue) under doubling of a genus--$2$ pencil with $2$ different pairs of base points, marked with $\bf{\times}$ in the figure.} 
     \label{figure0}
\end{figure}

Topologically, the doubling procedure can be interpreted roughly as follows: First break up the symplectic pencil $(X,f)$ into two pieces; the regular neighborhood of a smooth fiber (the ``concave piece''), which is a symplectic disk bundle of degree $m$ over a genus--$g$ surface, and its complement (the ``convex piece''). Describe a standard higher genus pencil on the former piece, based on $g$ and $m$ ---the model for which comes from holomorphic pencils. Then extend it over the convex piece, which now becomes a subpiece of the new pencil $(X,f')$. By the enterprise of \cite{AurouxKatzarkov}, the monodromy of $f'$ is explicitly determined by that of $f$ \cite{AurouxKatzarkov}[Theorem~4]. In particular, we observe  the following: any \emph{separating} Dehn twists in the positive factorization of $f'$ comes, in a one-to-one fashion, from  a separating Dehn twist in the positive factorization of $f$, which bounds a subsurface on $F$ that does not contain a base point of $f$. (See Figure~\ref{figure0} for an example.)


\section{Constructions} 

It is time to spell out our recipe in full detail. Take a positive factorization for a genus $g \geq 2$ pencil $(X,f)$ with base locus $B$ and regular fiber $F$, whose separating vanishing cycles are \,$a_1, \ldots, a_r$. Let $B'=\{b'_j\}$ and $B''=\{b''_j\}$ be two subsets of the base locus $B$ of $f$ (where the base points correspond to boundary components in the positive factorization), so that $|B'|=|B''| \leq 2g-2$. Set $K'_i, K''_i \in \{0,1\}$ to be the number of components in  $F \setminus a_i$  which do not intersect $B'$ (resp. $B''$).  \linebreak If $K'=\sum_{i=1}^{r'} K'_i$ is unequal to $K''=\sum_{i=1}^{r''} K''_i$, we are game. 

We blow-up all the base points of $f$ in $B \setminus B'$ (resp.\,$B \setminus B''$)  and double the resulting pencil on \,\mbox{$\widetilde{X}=X \, \# \,  |B \setminus B'| \, \CPb$} with $m=|B'|$ base points, which naturally come from the inclusion of $X$ minus the blown-up base points into $\widetilde{X}$. In this way, we arrive at a new pencil $(\widetilde{X},f')$ (resp.\,$(\widetilde{X}, f'')$). As we reviewed earlier, each one of the pencils $f'$ and $f''$ on $\widetilde{X}$  has genus $2g+m-1$ and $4m$ base points, whereas the number of reducible fibers they have are $K'$ and $K''$, respectively. (Figure~1 illustrates $2$ distinct partial doublings of sorts.) For each $0 \leq m_0 < m$, we get a pair of inequivalent pencils by blowing-up $m_0$ base points of both pencils, and for $m_0=m$, we get a pair of inequivalent fibrations. 

\smallskip

It remains to find the main ingredient for our recipe; a pencil $(X,f)$ as prescribed above, on (blow-ups of) rational and ruled surfaces.

\smallskip
\subsection{Inequivalent pencils and fibrations on rational surfaces} \

\begin{figure}[ht!]
 \centering
     \includegraphics[width=12.5cm]{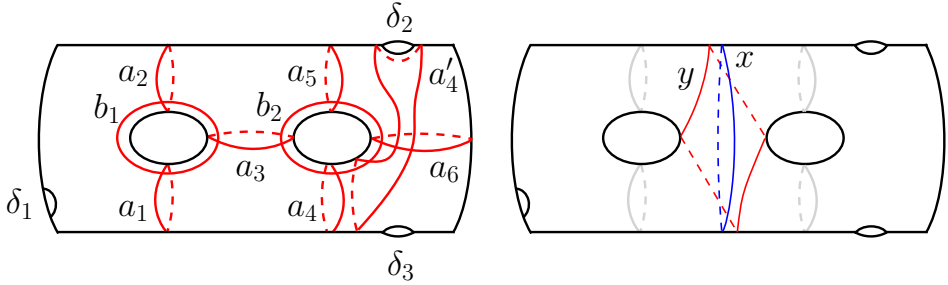}
     \caption{Twist curves $a_i, b_j$ and $a'_4$ on the surface $\Sigma_{2}^3$. The boundary components are $\delta_1, \delta_2, \delta_3$. The new curves $x,y$ are given in the second copy of $\Sigma_2^3$, where the lantern sphere sits in the middle. The separating curve $x$ is in blue.}
     \label{figure1}
\end{figure}

There is the following positive factorization of the boundary multi-twist in $\Gamma_2^3$, which was found by Sinem Onaran in \cite{Onaran}:
\begin{equation*}
\begin{array}{ccl}
t_{\delta_1} t_{\delta_2} t_{\delta_3} &=&  t_{a_3} t_{b_1} (t_{a_1}t_{a_2}t_{a_3}t_{b_1})^2 t_{a_3}  t_{\psi(b_2)} t_{a_3} t_{b_2} t_{a_4} t_{a_5} t_{a_3} t_{b_2} t_{a'_4} t_{a_6} \, .
\end{array}
\end{equation*}
Here the twist curves $a_i, b_j$ and $a'_4$ are as in Figure~\ref{figure1}, and $\psi = t_{a_5}^{-1} t_{a_4}^{-1}$. \linebreak (See \cite{Onaran}[Figure~6]. This is a lift of a well-known positive factorization of identity in $\Gamma_2$ \cite{BirmanHilden}, which corresponds to a hyperelliptic Lefschetz fibration on $\CP \# 13 \CPb$.)

We rewrite the above expression as
\begin{equation*}
\begin{array}{ccl}
t_{\delta_1} t_{\delta_2} t_{\delta_3} &=&
t_{a_3} t_{b_1} t_{a_1}t_{a_2}t_{a_3} t_{b_1} t_{a_1} t_{a_2}  t_{a_4} t_{a_5}  
\psi (t_{a_3}t_{b_1} t_{a_3} t_{\psi(b_2)} t_{a_3} t_{b_2}) \psi^{-1}
t_{a_3} t_{b_2} t_{a'_4} t_{a_6}  \\ &=&
t_{a_3} t_{b_1} t_{a_1}t_{a_2}t_{a_3} t_{b_1} \underline{t_{a_1} t_{a_2}  t_{a_4} t_{a_5} } 
t_{\psi(a_3)} t_{\psi(b_1)} t_{\psi(a_3)} t_{\psi^2(b_2)} t_{\psi(a_3)} t_{\psi(b_2)} 
t_{a_3} t_{b_2} t_{a'_4} t_{a_6}  
\end{array}
\end{equation*}
and then apply a lantern substitution along the underlined subword to obtain a new positive factorization:
\begin{equation*}
\begin{array}{ccl}
t_{\delta_1} t_{\delta_2} t_{\delta_3}  &=&
t_{a_3} t_{b_1} t_{a_1}t_{a_2}t_{a_3} t_{b_1} \underline{t_{a_3} t_{x}  t_{y}} 
t_{\psi(a_3)} t_{\psi(b_1)} t_{\psi(a_3)} t_{\psi^2(b_2)} t_{\psi(a_3)} t_{\psi(b_2)} 
t_{a_3} t_{b_2} t_{a'_4} t_{a_6} \, ,
\end{array}
\end{equation*}
which involves new twist curves $x, y$ given in Figure~\ref{figure1}. Note that $x$ is the only separating twist curve here.  Importantly, the ``lantern sphere'' we utilized did not contain the base points. The final positive factorization we obtained above describes a genus--$2$ Lefschetz pencil $(X,f)$ with $3$ base points.

The Seiberg-Witten adjunction inequality implies that in any closed symplectic $4$--manifold, except for rational or ruled surfaces, we have \,$2g(F)-2=-\eu(F) \geq [F]^2$ for any connected symplectic surface $F$ \cite{Ta, LiLiu}. However, $(X,f)$, equipped with the Gompf-Thurston symplectic form, has a symplectic fiber $F$ of genus $g(F)=2$ and self-intersection $[F]^2=3$ (the number of base points), violating the inequality. So we conclude that $X$ is a rational or ruled surface. We can then determine its exact diffeomorphism type through a simple calculation of its algebraic topological invariants. The first integral homology group $H_1(X)$, which is isomorphic to the quotient of $H_1(F)$ by the subgroup normally generated by  homology classes of vanishing cycles, is easily seen to be trivial: the vanishing cycles $a_1, b_1, a_3$ and $b_2$ already kill the whole group. Since the Euler characteristic of $X$ is given by \,$\eu(X)=4-4g+|C|-|B|$, where $C$ and $B$ are the critical locus and the base locus, we calculate $\eu(X)=12$. It follows that $X \cong \CP \# 9 \CPb$. 

Now let $b_j$ be the base point of $(X,f)$ corresponding to the boundary components $\delta_j$ in the factorization, for $j=1,2,3$. For $B'=\{b_2, b_3\}$ and $B''=\{b_1, b_2\}$, our recipe generates a pair of inequivalent genus--$5$ pencils \mbox{$f'$ and $f''$} on the rational ruled surface \mbox{$\widetilde{X} \cong \CP \# 10 \,\CPb \cong S^2\x S^2 \#\,9 \,\CPb$} with $8$ base points, where $f'$ has $1$ reducible fiber but $f''$ has none. (As illustrated in Figure~\ref{figure0}.)

\smallskip
\subsection{Inequivalent pencils and fibrations on irrational ruled surfaces} \

\begin{figure}[ht!]
 \centering
     \includegraphics[width=13cm]{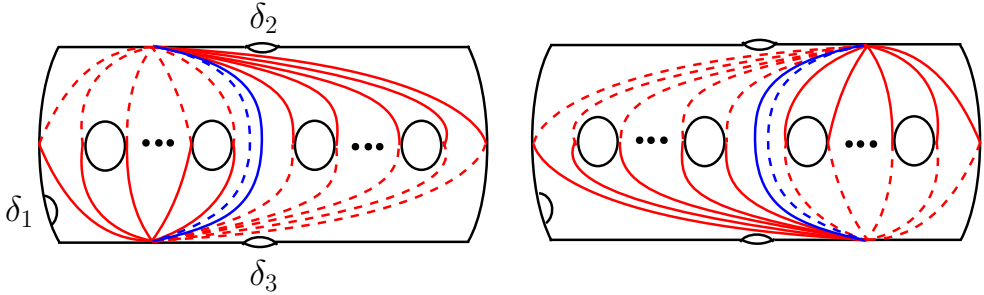}
     \caption{From left to right: twist curves $B_{0,1}, B_{1,1}, \ldots, B_{h,1}, C_1$, and $B_{0,2}, B_{1,2}, \ldots, B_{h,2}, C_2$ on (identical copies of) the surface $\Sigma_{2h}^3$. The boundary components are $\delta_1, \delta_2, \delta_3$. The separating curves $C_1$ and $C_2$ are given in blue.}
     \label{figure2}
\end{figure}

We next consider the following positive factorization of the boundary multi-twist in $\Gamma_{2h}^3$, for each $h \geq 1$, due to Noriyuki Hamada \cite{Hamada}:
\begin{equation*}
\begin{array}{ccl}
t_{\delta_1} t_{\delta_2} t_{\delta_3} &=& t_{B_{0,1}} t_{B_{1,1}} \cdots t_{B_{h,1}} t_{C_1} t_{B_{0,2}} t_{B_{1,2}} \cdots t_{B_{h,2}} t_{C_2} \, ,
\end{array}
\end{equation*}
where the curves $B_{i,k}, C_k$ for $i=0,1, \ldots, h$ and $k=1,2$ are as in Figure~\ref{figure2}. 
(This is the positive factorization ``$W_{IIA}$ for even genus'' in \cite{Hamada}[Figure~23] with only $3$ boundary components we picked and relabeled here. It is a lift of a well-known positive factorization of identity in $\Gamma_2$ \cite{Matsumoto, Korkmaz, Cadavid}, which corresponds to a hyperelliptic Lefschetz fibration on $S^2 \x \Sigma_h \# 4 \CPb$.) Corresponding to this positive factorization is a genus--$2h$ pencil $(X,f)$ with $3$ base points, where \mbox{$X\cong S^2 \x \Sigma_h \# \, \CPb$.} Note that $f$ has $2$ separating vanishing cycles. 

Once again, let $b_j$ be the base point of $(X,f)$ corresponding to the boundary component $\delta_j$ in the factorization, for $j=1,2,3$. We see that for $B'=\{b_2, b_3\}$ and $B''=\{b_1, b_2\}$, our recipe generates a pair of inequivalent genus--$(4h+1)$ pencils \mbox{$f'$ and $f''$} on the irrational ruled surface $\widetilde{X} \cong S^2 \x \Sigma_h \# \, 2\, \CPb \cong S^2 \widetilde{\x} \Sigma_h \# \, 2\, \CPb$ with $8$ base points, where $f'$ has $2$ reducible fibers but $f''$ has $1$.

\smallskip
\subsection{Generalizations.} \ 

We end with a few remarks on our recipe for constructing inequivalent pencils and fibrations in general:

1. From each pair of inequivalent pencils $f'$ and $f''$ we obtained, we can easily derive infinitely many more pairs of inequivalent pencils and fibrations of arbitrarily high genera. This is  achieved by taking further doubles of  $f'$ and $f''$; the base points of the new pencil will lie only on one component of each reducible fiber in $f'$ or $f''$ (see the monodromy in \cite{AurouxKatzarkov}[Theorem~4]), so the number of reducible fibers will remain different. 

 2. One can produce more than just a pair of inequivalent pencils and fibrations, whenever there is  an input pencil $(X,f)$ with several subcollections of base points $B', B'', \ldots, B^{(n)}$ of equal rank, which have different numbers $K', K'', \ldots, K^{(n)}$ of reducible fiber components they miss. For instance, a further lift of Hamada's positive factorization to $\Gamma_{2h}^4$ in \cite{Hamada} can be used to generate yet another pencil $f'''$ with no reducible fibers. 

 3. Our recipe can be improved in several ways. For example, we only compared the \emph{number }of reducible fibers, yet one can as well compare the \emph{topological types} of reducible fibers (determined by the genera of the subsurfaces they bound). It seems plausible that one can produce \emph{arbitrarily many} inequivalent fibrations and pencils by employing doubling \emph{sequences} as in \cite{BaykurJSG} for pencils on blow-ups of the same manifold, which contain  different topological types of reducible fibers.

\vspace{0.1in}
\noindent \textit{Acknowledgments.} These examples of inequivalent pencils and fibrations  came to life during my visit to Seoul, Korea in 2015, and were first announced during my lectures at the  highly stimulating \emph{2015 Seoul National University Topology Winter School}. This work is supported by the NSF Grant DMS-$1510395$.

\end{document}